\documentclass[a4paper,11pt]{article}
\usepackage{hyperref}
\hypersetup{
colorlinks=true,
linktoc=all,
linkcolor=black,
citecolor=blue,
}
\usepackage[draft]{changes}
\usepackage{cite}
\usepackage{setspace}
\usepackage{indentfirst}
\usepackage{amssymb,amsthm}
\usepackage{amsmath}
\allowdisplaybreaks
\newenvironment{pf}{\noindent {\bf Proof.}}{\hfill\qedsymbol\par\medskip}
\newenvironment{pf1}{\noindent {\bf Proof of Theorem \ref{mainresult}.}}{\hfill\qedsymbol\par\medskip}

\newtheorem{theorem}{Theorem}[section]
\newtheorem{lemma}[theorem]{Lemma}

\begin{document}
\begin{spacing}{1}

\title{Ramsey goodness of large generalized fans}
\date{}

\author{
Shaonan Mi\footnote{
College of Mathematical Sciences, Harbin Engineering University, Harbin 150001, China.
Email: {\tt mishaonan@hrbeu.edu.cn}.}
\;\;\;\; Ye Wang\footnote{
Corresponding author.
College of Mathematical Sciences, Harbin Engineering University, Harbin 150001, China.
Email: {\tt ywang@hrbeu.edu.cn}.
Supported in part by NSFC (No.\ 12471323).
}
}
\maketitle

\begin{abstract}
Given graphs $G$ and $H$, the Ramsey number $r(G,H)$ is the
smallest integer $N$ such that any red-blue edge coloring of
$K_N$ contains a red copy of $G$ or a blue copy of $H$. If $H$ is
connected, Burr's lower bound gives
$r(G,H)\geq (\chi(G)-1)(|V(H)|-1)+s(G),$
where $s(G)$ is the minimum size of a color class over all proper
$\chi(G)$-colorings of $G$. The graph $H$ is called $G$-good if equality
holds in this bound. Liu and Li [J. Graph Theory, 2025] characterized the graphs $G$ with $s(G)=1$ for which
$K_p+nH$ is $G$-good for $p\geq 2$ and sufficiently large $n$, and posed the corresponding problem
for $p=1$. We solve this problem by showing that 
for a graph $G$ with $\chi(G)=k+1\geq 2$ and $s(G)=1$, a graph $H$ of order $h$ and sufficiently large $n$, 
$K_1+nH$ is $G$-good if and only if
\[G\subseteq
\begin{cases}
 mK_2+K_{k-1}(m), & \text{if $nh$ is odd},\\[2mm]
 (mK_2\cup K_{1,2})+K_{k-1}(m),
    & \text{if $nh$ is even},
\end{cases}
\]
for some positive integer $m$.
\end{abstract}

\section{Introduction}

Given graphs $G$ and $H$, the Ramsey number $r(G,H)$ is the smallest integer $N$ such that every red-blue edge coloring of $K_N$ contains a red copy of $G$ or a blue copy of $H$.
Let $\chi(G)$ denote the chromatic number of $G$, and define the chromatic surplus $s(G)$ to be the minimum size of a color class among all proper $\chi(G)$-colorings of $G$.
For vertex-disjoint graphs $G$ and $H$, we write $G\cup H$ and $G+H$ for the union and the join of $G$ and $H$, respectively.
Denote by $nG$ the union of $n$ vertex-disjoint copies of $G$.
Burr \cite{Bur81} proved that every connected graph $H$ with $|V(H)|\geq s(G)$ satisfies
\begin{equation}\label{burr}
r(G,H)\geq(\chi(G)-1)(|V(H)|-1)+s(G).
\end{equation}
Following Burr and Erd\H{o}s \cite{BE83}, we say that $H$ is $G$-good if equality holds in \eqref{burr}.

Chv\'atal's theorem that every tree is $K_p$-good is a fundamental example of Ramsey goodness \cite{Chv77}. 
Burr and Erd\H{o}s initiated the systematic study of this subject \cite{BE83}, 
followed by results for various graph families and general structural methods \cite{BFRS83,NR09}.

A particularly relevant line of work concerns Ramsey goodness for generalized fans.
Li and Rousseau \cite{LR96} proved that, for any fixed graphs $F$ and $H$, $K_1+nH$ is $(K_2+F)$-good for all sufficiently large $n$.
Lin, Li and Dong \cite{LLD10}, and later Chung and Lin \cite{CL25}, 
determined goodness for several classes of complete multipartite graphs, revealing essential parity phenomena. 
Jiang and Hou \cite{JH25} proved that $K_1+nK_t$ is $G$-good for
every edge-critical graph $G$ with $\chi(G)\geq3$ and all sufficiently
large $n$. Zhang and Chen \cite{ZC26} obtained improved sufficient
ranges for $K_1+nH$ to be $K_p$-good.

We are primarily concerned with structural characterizations of the graphs $G$ for which a given family is $G$-good.
For positive integers $t_1,\ldots,t_p$ with $p\ge 1$, let $K_{t_1,\ldots,t_p}$ denote the complete $p$-partite graph whose parts have sizes $t_1,\ldots,t_p$, respectively.
When all parts have the same size $t$, we denote this graph by $K_p(t)$.
Burr and Faudree \cite{BF93} obtained the following structural characterization for large trees.

\begin{theorem}[Burr and Faudree~\cite{BF93}]
Let $G$ be a fixed graph with $\chi(G)=k+1\ge2$ and
$s(G)=1$.
Then all sufficiently large trees are $G$-good if and
only if
\[
G\subseteq mK_2+K_{k-1}(m)
\]
for some positive integer $m$.
\end{theorem}
Liu and Li \cite{LL25} showed that the same structural condition characterizes Ramsey goodness for large generalized books.
\begin{theorem}[Liu and Li~\cite{LL25}]
Let $G$ be a fixed graph with $\chi(G)=k+1\ge2$ and
$s(G)=1$, let $H$ be a fixed graph of order $h$,
and let $p\ge2$ be a fixed integer.
Then large $K_p+nH$ are $G$-good
if and only if
\[
G\subseteq mK_2+K_{k-1}(m)
\]
for some positive integer $m$.
\end{theorem}
Liu and Li~\cite{LL25} also posed the corresponding
characterization problem for $p=1$.
In this case, they proved that
$G\subseteq mK_2+K_{k-1}(m)$
for some positive integer $m$ is sufficient for
$K_1+nH$ to be $G$-good for all sufficiently large $n$,
whereas
$G\subseteq (mK_2\cup K_{1,2})+K_{k-1}(m)$
for some positive integer $m$ is necessary.
In particular, Chung and Lin~\cite{CL25} proved that
$K_1+nH$ is $(K_1+C_4)$-good for all sufficiently large
$n$ with $n|V(H)|$ even.
In this paper, we solve the characterization
problem by distinguishing the two possible
parities of $n|V(H)|$.

\begin{theorem}\label{mainresult}
For a graph $G$ with $\chi(G)=k+1\geq 2$ and $s(G)=1$, a graph $H$ of order $h$ and sufficiently large $n$, 
$K_1+nH$ is $G$-good if and only if
\[
G\subseteq
\begin{cases}
mK_2+K_{k-1}(m),
& \text{if $nh$ is odd},\\[1mm]
(mK_2\cup K_{1,2})+K_{k-1}(m),
& \text{if $nh$ is even},
\end{cases}
\]
for some positive integer $m$.
\end{theorem}

Theorem \ref{mainresult} resolves the case $s(G)=1$ of the characterization
problem posed by Zhang and Chen~\cite[Problem~1.6]{ZC26}.
The structural characterizations of Burr and
Faudree~\cite{BF93} for large trees and of Liu and
Li~\cite{LL25} for large generalized books, as well as our result,
all assume that $s(G)=1$.
It is natural to seek analogous structural characterizations
beyond chromatic surplus one, starting with connected
graphs $G$ satisfying $s(G)=2$.
Even this first case exhibits arithmetic restrictions
beyond parity.
Let $P_\ell$ denote the path on $\ell$ vertices and let
$F_n=K_1+nK_2$.
As a special case of the results of Li and Ning~\cite{LiNing},
for every $n\ge3$,
\[
r(P_4,F_n)=
\begin{cases}
2n+2, & n\not\equiv2\pmod3,\\
2n+3, & n\equiv2\pmod3.
\end{cases}
\]
Since $\chi(P_4)=s(P_4)=2$, for $n\ge 3$,  $F_n$ is $P_4$-good
if and only if $n\not\equiv2\pmod3$.
Thus, in contrast to the case $s(G)=1$ in Theorem~\ref{mainresult},
goodness need not stabilize for large $n$ even when the
parity of $n|V(H)|$ is fixed.
It is therefore natural to ask how these arithmetic conditions interact with the structure of \(G\), 
and whether the characterization can be extended to the case \(s(G)>1\).

\medskip

\section{Proof of the main result}

Before the proof, we list and prove some preliminary results.
To state these results, we first introduce some standard notation.
For $v\in V(G)$ and $S\subseteq V(G)$, 
let $N_G(v)$ be the neighborhood of $v$ in $G$, 
write $d_G(v)=|N_G(v)|$ and $d_G(v,S)=|N_G(v)\cap S|$,
and let $G[S]$ denote the subgraph of $G$ induced by $S$.
Denote the minimum and maximum degrees of $G$ by $\delta(G)$ and $\Delta(G)$, respectively.
We shall use the following structural lemma of Chung and Lin \cite[Lemma 4.4]{CL25}.

\begin{lemma}[Chung and Lin \cite{CL25}]\label{CLstructure}
For any fixed graph $H$ and integers $p\geq2$ and $b\ge 1$,
there exists $\delta>0$ such that the following holds for all $n\geq pb^2/\delta$.
Let $a_1,\ldots,a_p$ be positive integers with
$1=a_1\leq a_2\leq\cdots\leq a_{p-1}\leq b$ and $a_p\leq\delta n$.
If $\Gamma$ is a $K_{a_1,\ldots,a_p}$-free graph on at least $(p-1)n|V(H)|$ vertices and $\overline{\Gamma}$ is $(K_1+nH)$-free, then $V(\Gamma)$ has a partition $V_1\cup\cdots\cup V_{p-1}$ such that $\Delta(\Gamma[V_i])\leq a_2-1$, for $i=1,2,\ldots,p-1$.
\end{lemma}

We also record the following elementary estimate,
which follows from the argument of Burr, Erd\H{o}s and Spencer \cite[Lemma~1]{BES75}.

\begin{lemma}[Burr, Erd\H{o}s and Spencer \cite{BES75}]\label{packingbound}
For any graphs $G$ and $H$ and any positive integer $n$,
$r(G,nH)\leq(n-1)|V(H)|+r(G,H).$
\end{lemma}

The equitable-coloring form of the Hajnal-Szemer\'edi theorem says that, for every positive integer $r$, every graph of maximum degree at most $r$ admits a proper equitable $(r+1)$-coloring, in which the sizes of any two color classes differ by at most one; see Hajnal and Szemer\'edi \cite{HS70}, and Kierstead and Kostochka \cite{KK08} for a later short proof.
Applying the Hajnal-Szemer\'edi theorem to the complement gives the following lemma.

\begin{lemma}[Hajnal and Szemer\'edi \cite{HS70}, and Kierstead and Kostochka \cite{KK08}]\label{cliquefactor}
For $h\geq1$ and $n\geq2$,
if $\Gamma$ is a graph of order $nh$ and $\delta(\Gamma)\geq nh-n$, then $\Gamma$ contains a $K_h$-factor, that is, $n$ vertex-disjoint copies of $K_h$ covering $V(\Gamma)$.
\end{lemma}


The following result extracts a parity-dependent fixed subgraph from a bounded-degree graph.

\begin{lemma}\label{paritydefect}
Let $H$ be a graph of order $h$, and let $t,d$ be positive integers.
There exists a positive integer $n_1=n_1(H,t,d)$ such that, for every
integer $n\ge n_1$, the following holds.
If $F$ is a graph of order at least $nh+1$ with
$\delta(F)\ge 1$ and $\Delta(F)\le d$, and $\overline F$ is
$(K_1+nH)$-free, then $tK_2\subseteq F$ if $nh$ is odd, and
$tK_2\cup K_{1,2}\subseteq F$ if $nh$ is even.
\end{lemma}

\begin{pf}
Let
\begin{equation}\label{defectthreshold}
n_1=\max\Big\{d+1, \Big\lceil \frac{(2t+1)d}{h}\Big\rceil \Big\}.
\end{equation}

Suppose first that $nh$ is odd. Let $M$ be a maximum matching in $F$.
Suppose to the contrary that $M$ has at most $t-1$ edges, and then consider the vertices not covered by $M$.
Since $M$ is maximal and $\delta(F)\geq 1$, 
every vertex in $V(F)\setminus V(M)$ is adjacent to at least one vertex in $V(M)$.
Hence $|V(F)\setminus V(M)|\leq (d-1)|V(M)|$, and therefore $|V(F)|\leq d|V(M)|\leq2(t-1)d$, contrary to \eqref{defectthreshold} and $|V(F)|\geq nh+1$.
Thus $tK_2\subseteq F$.

Now suppose that $nh$ is even.
We first show that $F$ contains a copy of $K_{1,2}$.
If not, since $\delta(F)\ge 1$,  each vertex in $F$ has only one neighbor, which implies that $F$ is a perfect matching.
As $nh$ is even and $|V(F)|\geq nh+1$, we have $|V(F)|\geq nh+2$.
For a vertex $x\in V(F)$, $d_{\overline F}(x)=|V(F)|-2\geq nh$.
Choose a set $W\subseteq N_{\overline F}(x)$ with $|W|=nh$.
For every $w\in W$,  
\[
d_{\overline F[W]}(w)=nh-1-d_{F[W]}(w)\geq nh-1-d\geq nh-n.
\]
Lemma \ref{cliquefactor} yields a $K_h$-factor in $\overline F[W]$.
Each of its $n$ cliques contains a copy of $H$, and these copies together with the center $x$ form $K_1+nH$ in $\overline F$, a contradiction.
So $F$ contains a copy of $K_{1,2}$.

Let $U$ be the vertex set of a copy of $K_{1,2}$ contained in $F$.
It suffices to show that $V(F)\setminus U$ induces a subgraph containing a matching with $t$ edges.
Let $M$ be a maximum matching in the subgraph induced by $V(F)\setminus U$.
Suppose to the contrary that $M$ has at most $t-1$ edges, and then consider the vertices in $V(F)\setminus (V(M)\cup U)$.
Since $M$ is maximal and $\delta(F)\geq 1$, 
every vertex in $V(F)\setminus (V(M)\cup U)$ is adjacent to at least one vertex in $V(M)\cup U$.
Hence 
$|V(F)\setminus (V(M)\cup U)|\leq (d-1)(|V(M)|+|U|)$. 
Equivalently, 
$|V(F)|\leq d(|V(M)|+|U|)\le (2t+1)d$, 
contrary to \eqref{defectthreshold}.
So $tK_2\cup K_{1,2}\subseteq F$.
\end{pf}

\begin{pf1}
We prove the necessity and sufficiency separately.

\medskip
\noindent\textbf{($\Rightarrow$) Necessity.}
Suppose that $r(G,K_1+nH)=knh+1$.
If $nh$ is odd,
consider a red-blue edge coloring of $K_{knh+1}$ such that the red graph
\[
R=\frac{nh+1}{2}K_2+K_{k-1}(nh).
\]
Then the blue complement is the union of $k-1$ copies of $K_{nh}$ and the complement of $\frac{nh+1}{2}K_2$.
Every vertex in every blue component has degree $nh-1$, whereas the center of $K_1+nH$ has degree $nh$.
Hence the blue graph contains no copy of $K_1+nH$.
The assumed Ramsey equality forces a red copy of $G$ in $R$.
Since $G$ is fixed, $G\subseteq mK_2+K_{k-1}(m)$ for some positive integer $m$.

If $nh$ is even,
consider a red-blue edge coloring of $K_{knh+1}$ such that the red graph
\[
R=\Big(\frac{nh-2}{2}K_2\cup K_{1,2}\Big)+K_{k-1}(nh).
\]
In the blue component, every vertex has degree at most $nh-1$, 
and hence the blue graph contains no $K_1+nH$.
Thus, $G\subseteq(mK_2\cup K_{1,2})+K_{k-1}(m)$.

\medskip
\noindent\textbf{($\Leftarrow$) Sufficiency.}
Since $s(G)=1$, fix a proper $(k+1)$-coloring of $G$ whose color-class sizes are $1=a_1\leq a_2\leq\cdots\leq a_{k+1}$.
Let $Q=K_{a_1,a_2,\ldots,a_{k+1}}$ and $\ell=r(Q,H)$.
Note that $G\subseteq Q$.

Let $N=knh+1$, and consider a red-blue edge coloring of $K_N$ with red and blue spanning graphs $R$ and $B$, respectively.
Suppose, for a contradiction, that $R$ contains no copy of $G$ and $B$ contains no copy of $K_1+nH$.
Since $G\subseteq Q$, the red graph $R$ is also $Q$-free.

We first bound the blue degrees.
If some vertex $v$ satisfies $d_B(v)\geq(n-1)h+\ell$, choose a set $W\subseteq N_B(v)$ with  $|W|=(n-1)h+\ell$.
By Lemma \ref{packingbound}, $W$ induces a red copy of $Q$ or a blue copy of $nH$,
which yields a red copy of $Q$ or a blue copy of $K_1+nH$ along with $v$, a contradiction.
So for any $v\in V(B)$,
\begin{equation}\label{bluedegree}
d_B(v)\leq(n-1)h+\ell-1.
\end{equation}

Note that $R$ is $Q$-free, $B$ is $(K_1+nH)$-free, and $|V(R)|=N=knh+1\geq kn|V(H)|$.
By Lemma \ref{CLstructure},
there exists $\delta>0$ such that for all $n\geq (k+1)a_k^2/\delta$,
$V(R)$ has a partition $V_1\cup\cdots\cup V_k$ such that
\begin{equation}\label{internaldegree}
\Delta(R[V_i])\leq a_2-1,
\qquad\text{for }1\le i\le k.
\end{equation}
Let $c_0=\ell+a_2-h-1$.
For any $v\in V_i$, equations \eqref{bluedegree} and \eqref{internaldegree} give
\[
|V_i|-a_2
\leq d_B(v,V_i)
\leq d_B(v)
\leq nh+\ell-h-1,
\]
which yields $|V_i|\leq nh+c_0$.
Using this upper bound for the other $k-1$ parts, we obtain
\begin{equation}\label{partlower}
|V_i|
\geq N-(k-1)(nh+c_0)
=nh+1-(k-1)c_0,
\qquad\text{for }1\le i\le k.
\end{equation}
If $i\neq j$ and $v\in V_i$, then
\begin{equation}\label{crossblue}
\begin{aligned}
d_B(v,V_j)
&\leq d_B(v)-d_B(v,V_i)\\
&\leq nh+\ell-h-1-(|V_i|-a_2)\\
&\leq \ell-h+a_2-2+(k-1)c_0\\
&\leq kc_0.
\end{aligned}
\end{equation}
Thus \eqref{crossblue} holds for all distinct $i,j\in[k]$ and every $v\in V_i$.

Let $S$ be a largest part of the partition $V_1\cup\cdots\cup V_k$.
Since $N=knh+1$, we have $|S|\geq\lceil N/k\rceil=nh+1$.
If $a_2=1$, then \eqref{internaldegree} shows that $R[S]$ is edgeless,
and the blue graph $B[S]$ is a clique on at least $nh+1$ vertices and contains a copy of $K_1+nH$, a contradiction.
Now assume that $a_2\geq2$.

We claim that $\delta(R[S])\geq1$.
Otherwise, let $x$ be an isolated vertex of $R[S]$ and choose a set $W\subseteq S\setminus\{x\}$ with $|W|=nh$.
Every edge from $x$ to $W$ is blue, and \eqref{internaldegree} gives $\delta(B[W])\geq nh-a_2\geq nh-n$.
Lemma \ref{cliquefactor} gives a blue $K_h$-factor in $W$.
Together with the center $x$, we obtain a blue copy of $K_1+nH$, a contradiction.
Thus $1\leq\delta(R[S])\leq\Delta(R[S])\leq a_2-1$.

Since $R[S]$ is a graph of order at least $nh+1$ and  $B[S]$ is $(K_1+nH)$-free,
by Lemma \ref{paritydefect}, for sufficiently large $n$,
it gives 
\begin{equation}\label{extracteddefect}
\begin{cases}
mK_2\subseteq R[S],&\text{if $nh$ is odd},\\
mK_2\cup K_{1,2}\subseteq R[S],&\text{if $nh$ is even},
\end{cases}
\end{equation}
for any positive integer $m$.
This supplies the parity-dependent part of the corresponding graph in Theorem \ref{mainresult},
and the remaining $k-1$ parts will be embedded next.

Relabel the partition so that $S=V_1$.
The assumed embedding of $G$ induces a vertex partition $V(G)=A\cup C_1\cup\cdots\cup C_{k-1}$ such that
\[
\begin{cases}
G[A]\subseteq mK_2,&\text{if $nh$ is odd},\\
G[A]\subseteq mK_2\cup K_{1,2},&\text{if $nh$ is even},
\end{cases}
\]
and every $C_j$ is independent in $G$ with $|C_j|\leq m$.
Use the corresponding red copy in \eqref{extracteddefect} to embed $A$ into $V_1$.

We embed $C_j$ into $V_{j+1}$ successively for $1\le j\le k-1$.
Suppose that $A\cup C_1\cup \ldots \cup C_{j-1}$ have already been embedded, 
and let $X_j$ be the set of the images in $V_1\cup V_2\cup\ldots\cup V_{j}$ for $2\le j\le k-1$.
If $j=1$, set $X_1$ to be the set of the images of $A$ in $V_1$.
Every vertex of $X_j$ lies outside $V_{j+1}$, and $|X_j|\leq km+3$.
By \eqref{crossblue}, each $x\in X_j$ excludes at most $kc_0$ vertices of $V_{j+1}$ from the common red neighborhood.
By \eqref{partlower} and the union bound,
\[
\begin{aligned}
\Big| \bigcap_{x \in X_j} \bigl( N_R(x) \cap V_{j+1} \bigr) \Big|
&\geq |V_{j+1}|-\sum_{x\in X_j}d_B(x,V_{j+1})\\
&\geq nh+1-(k-1)c_0-k|X_j|c_0\\
&\geq nh+1-(k-1)c_0-k(km+3)c_0\\
&\ge m,
\end{aligned}
\]
if $n$ is chosen sufficiently large.
We can therefore choose $|C_j|\leq m$ distinct vertices from this common red neighborhood in $V_{j+1}$.
Repeating the procedure embeds every $C_j$.

Note that all required edges inside $A$ are supplied by the corresponding red copy in \eqref{extracteddefect}, all required edges between distinct sets are red by construction, and no edge is required inside any $C_j$.
So we obtain a red copy of $G$, contrary to the assumption.
Thus $r(G,K_1+nH)\leq knh+1$ for sufficiently large $n$.
Together with the lower bound from \eqref{burr}, this completes the sufficiency argument.
\end{pf1}

\section*{Declaration on the use of AI}

During the preparation of this work, the authors used ChatGPT (OpenAI) to assist with mathematical exploration. Using Lemma 2.1, ChatGPT derived a result similar to Lemma 2.4 and showed that
for every fixed bipartite graph $G$ and integers $k\ge 2$ and $m\ge 3$, $F_n$ is $(G+K_1+K_{k-2}(m))$-good for all sufficiently large $n$ if and only if
$G\subseteq qK_2\cup K_{1,2}$ or $G\subseteq K_{2,q}$ for some positive integer $q$.
The authors subsequently refined and completed the arguments and established the necessary and sufficient conditions presented in this paper. The authors take full responsibility for the mathematical correctness and the final content of the manuscript.

\end{spacing}
\end{document}